\documentclass[10pt]{amsart}
\usepackage{amsfonts}
\usepackage{amssymb}
\usepackage{amsthm}
\begin{document}

%%%%%%%%%%%%%%%%%%%%%%%%%%%%%%%%%%%%
%  Macros                          %
%%%%%%%%%%%%%%%%%%%%%%%%%%%%%%%%%%%%
\newcommand{\Zthree}{\mathbb{Z}^3\,}
\newcommand{\Ztwo}{\mathbb{Z}^2\, }
\newcommand{\Zone}{\mathbb{Z}\, }
\newcommand{\Ttwo}{\mathbb{T}^2\, }
\newcommand{\Tone}{\mathbb{T}\, }
\newcommand{\Cone}{\mathbb{C}\,}
\newcommand{\Cn}{\mathbb{C}^n\,}
\newcommand{\Rone}{\mathbb{R}\,}
\newcommand{\Rtwo}{\mathbb{R}^2\,}
\newcommand{\Rn}{\mathbb{R}^n\,}
\newcommand{\nat}{\mathbb{N}\,}
\newcommand{\J}{\mathbb{J}\,}
\newcommand{\hs}{\mathcal{H}\,}             %Hilbert space%
\newcommand{\bh}{\mathcal{B}(\hs)\,}        %Bounded operators on H%
\newcommand{\f}{\mathcal{F}\,}
\newcommand{\fn}{\f^{-1}}
\renewcommand{\S}{\mathcal{S}_1\,}  %unit sphere%
\newcommand{\A}{\alpha}
\newcommand{\B}{\beta}
\newcommand{\D}{\delta}
\newcommand{\E}{\epsilon}
\renewcommand{\L}{\lambda}
\newcommand{\T}{\theta}
\renewcommand{\O}{\Omega}
\newcommand{\s}{\sigma}
\newcommand{\txt}[1]{\;\mbox{#1}\;}
\newcommand{\tx}[1]{\;\mbox{#1}\,}
\newcommand{\ip}[2]{\langle {#1},{#2} \rangle}
\newcommand{\nn}{\nonumber}
\newcommand{\dd}{\,\mathrm{d}}
\newcommand{\bld}[1]{\textbf{#1}}
\newcommand{\noi}{\noindent}
\newcommand{\bt}{\begin{theorem}}
\newcommand{\bl}{\begin{lemma}}
\newcommand{\bp}{\begin{proof}}
\newcommand{\bd}{\begin{definition}}
\newcommand{\bc}{\begin{corollary}}
\newcommand{\beq}{\begin{eqnarray}}
\newcommand{\beqn}{\begin{eqnarray*}}
\newcommand{\et}{\end{theorem}}
\newcommand{\el}{\end{lemma}}
\newcommand{\ep}{\end{proof}}
\newcommand{\ed}{\end{definition}}
\newcommand{\ec}{\end{corollary}}
\newcommand{\eeq}{\end{eqnarray}}
\newcommand{\eeqn}{\end{eqnarray*}}

%%%%%%%%%%%%%%%%%%%%%%%%%%%%%%%%%%%%
% More   Macros                    %
%%%%%%%%%%%%%%%%%%%%%%%%%%%%%%%%%%%%
\newcommand\Ec{{\mathcal{E}}}
\newcommand\Nats{{\mathbb N}}
\newcommand\Reals{{\mathbb R}}
\newcommand\ub{{\mathbf u}}
\newcommand\vb{{\mathbf v}}
\newcommand\Vc{{\mathcal{V}}}
\newcommand\xb{{\mathbf x}}
\newcommand\yb{{\mathbf y}}

%%%%%%%%%%%%%%%%%%%%%%%%%%%%%%%%%%%%%
%     Theorems                      %
%%%%%%%%%%%%%%%%%%%%%%%%%%%%%%%%%%%%%
\newtheorem{theorem}{Theorem}
\newtheorem{lemma}[theorem]{Lemma}
\newtheorem{proposition}[theorem]{Proposition}

\newtheorem*{nntheorem}{Theorem}
\newtheorem{corollary}[theorem]{Corollary}

\theoremstyle{remark}
\newtheorem{remark}[theorem]{Remark}
\theoremstyle{definition}
\newtheorem{definition}{Definition}

%%%%%%%%%%%%%%%%%%%%%%%%%%%%%%%%%%%%%%
%        Title                       %
%%%%%%%%%%%%%%%%%%%%%%%%%%%%%%%%%%%%%%

\title[Ellipsoidal Tight Frames]{Ellipsoidal Tight Frames and Projection
Decompositions of Operators}

\author[Dykema]{Ken Dykema}
\address{Department of Mathematics \\
 Texas
A\&M University \\ 3368 TAMU \\ College Station, TX 77843-3368}
\email{kdykema@math.tamu.edu}

\author[Freeman]{Dan Freeman}
\address{ Guilford College\\  5800 West Friendly Avenue \\
Greensboro, NC 27410}
\email{dfreeman@guilford.edu}

\author[Kornelson]{Keri Kornelson}
\address{Department of Mathematics \\
 Texas
A\&M University \\ 3368 TAMU \\ College Station, TX 77843-3368}
\email{keri@math.tamu.edu}

\author[Larson]{David Larson}
\address{Department of Mathematics \\
 Texas
A\&M University \\ 3368 TAMU \\ College Station, TX 77843-3368}
\email{larson@math.tamu.edu}

\author[Ordower]{Marc Ordower}
\address{Department of Mathematics \\ Randolph Macon Woman's
 College \\
2500 Rivermont Ave. \\
Lynchburg, VA
24503 }
\email{mordower@rmwc.edu}

\author[Weber]{Eric Weber}
\address{ Department of Mathematics, University of Wyoming,  P.O. Box
 3036
Laramie, WY 82071}
\email{esw@uwyo.edu}

\thanks{The research of the first, fourth, and sixth authors was supported in part by grants
from the NSF}


%%%%%%%%%%%%%%%%%%%%%%%%%%%%%%%%%%%%%%
%    Abstract                        %
%%%%%%%%%%%%%%%%%%%%%%%%%%%%%%%%%%%%%%
\begin{abstract}
  We prove the
existence of tight frames whose elements lie on an arbitrary
ellipsoidal surface within a real or complex separable Hilbert space $\hs$, and we analyze the set 
of attainable frame bounds.
In the case where $\hs$ is real and has finite dimension, we give an algorithmic
 proof. Our main tool in the infinite-dimensional
case is a result we have proven which concerns the
 decomposition of a positive
invertible operator into a strongly converging sum of (not necessarily mutually
orthogonal) self-adjoint projections.  This decomposition result implies the existence of tight
frames in the ellipsoidal surface  determined by the
positive operator.  In the real or complex finite dimensional case, this provides an alternate (but
 not
algorithmic) proof that every such surface contains tight frames with every prescribed length at
 least as large as
dim$\,\hs$. A corollary in both finite and infinite dimensions is that every positive
invertible  operator is the frame
operator for a spherical frame.
\end{abstract}

\maketitle

%%%%%%%%%%%%%%%%%%%%%%%%%%%%%%%%%%%%%%%%%%%%%%%%%%%%%%%%%%%%%%%%
%              Section                                         %
%%%%%%%%%%%%%%%%%%%%%%%%%%%%%%%%%%%%%%%%%%%%%%%%%%%%%%%%%%%%%%%%

\section*{Introduction}\label{intro}
Frames were first introduced by Dufflin and Schaeffer \cite{DS} in
1952 as a component in the development of non-harmonic Fourier
series, and a paper by Daubechies, Grossmann, and Meyer \cite{DGM}
in 1986 initiated the use of frame theory in signal processing. A
\textit{frame} on a separable Hilbert space $\hs$ is defined to be
a complete collection of vectors $\{x_i\}\subset \hs$ for which
there exists constants $0 < A \leq B$ such that for any $x \in
\hs$,
$$A\|x\|^2 \leq \sum_i |\ip{x}{x_i}|^2 \leq B \|x\|^2$$
The constants $A$ and $B$ are known as the \textit{frame bounds}. The collection is called a
 \textit{tight frame}
if $A=B$, and a \textit{Parseval frame} if $A=B=1$.  (In some of the existing literature, Parseval
 frames have
been called \textit{normalized tight frames}, however it should be noted that other authors have
 used the term
{\it normalized} to describe a frame consisting only of unit vectors.) The \textit{length} of a
 frame is the
number of vectors it contains, which cannot be less than the Hilbert space dimension.  References in
 the study of
frames include \cite{D}, \cite{HL}, and \cite{HW}.

Hilbert space frames are used in a variety of signal processing applications, often demanding
 additional
structure.  Tight frames may be constructed having specified length, components having a
 predetermined sequence of
norms, or with properties making them resilient to erasures. For examples, see \cite{BF},
 \cite{CKLT}, \cite{GKK}, and \cite{HP}.
 One area of rapidly advancing research lies in
describing tight frames in which all the vectors are of equal
norm, and thus are elements of a sphere. \cite{BF}, \cite{CKLT}

Since frame theory is geometric in nature, it is natural to ask which other surfaces in a finite or
 infinite
dimensional Hilbert space contain tight frames. By an {\it ellipsoidal surface} we mean the image of
 the unit
sphere $\S =\{x \,:\, \|x\|=1\}$ under a bounded invertible operator $T \in \bh$. Let $\Ec_T$ denote
 the
ellipsoidal surface $\Ec_T = T\S$.  A frame contained in $\Ec_T$ is called an {\it
 ellipsoidal
 frame}, and if it is tight it is called an ellipsoidal tight frame (ETF) for
 that surface.
We say that a frame bound K is {\it attainable} for $\Ec_T$
 if there is an ETF for $\Ec_T$ with frame bound K.  If an ellipsoid $\Ec$ is a
 sphere we will call a frame in $\Ec$  {\it spherical}.

Given an ellipsoid $\Ec$, we can assume $\Ec=\Ec_T$ where $T$ is a positive
 invertible operator.
Given $A$ an invertible operator, let $A^* = U|A^*|$ be the polar decomposition
 where $|A^*| = (AA^*)^{\frac{1}{2}}$.
Then $A = |A^*|\, U^*$.  By taking $T=|A^*|$, we  that $T\S = A\S$.  Moreover it
 is easily seen that the positive operator
$T$ for which $\Ec = \Ec_T$ is unique.

Throughout the paper, $\hs$ will be a separable real or complex Hilbert space and for $x,y,u \in
 \hs$, we will use
the notation $x \otimes y$ to denote the rank-one operator $u \mapsto \ip{u}{y} x$.  Note that
 $\|x\|=1$ implies
that $x \otimes x$ is a rank-1 projection.

Special thanks are given to our colleagues
Pete Casazza, Vern Paulsen (see Remark \ref{others}), and Nicolaas Spronk for useful conversations
concerning the material in this paper, and also to undergraduate
REU/VIGRE research students Emily King, Nate Strawn and Justin Turner for
taking part in discussions on ellipsoidal frames.  This research project began in an REU/VIGRE
 seminar course at Texas A\&M University in Summer 2002
 in which all the
co-authors were participants.

%%%%%%%%%%%%%%%%%%%%%%%%%%%%%%%%%%%%%%%%%%%%%%%%%%%%%%%%%%%%%%%%
%              Section                                         %
%%%%%%%%%%%%%%%%%%%%%%%%%%%%%%%%%%%%%%%%%%%%%%%%%%%%%%%%%%%%%%%%
\section {Theorems}\label{thms}
There are three theorems in this paper.  The first gives an elementary
construction of ETF's when $\hs = \Rone^n$, and is proved in Section
 \ref{rotate}.

%%%%%%%%%   Theorem 1   %%%%%%%%%%
 \bt\label{construct} Let $n,k\in\Nats$ with
$n\le k$, let $a_1,\ldots,a_n\ge0$ be such that $r:=\sum_1^na_j>0$
and consider the (possibly degenerate) ellipsoid
\[
\Ec=\{\xb=(x_1,\ldots,x_n)^t\in\Reals^n\mid\sum_1^na_jx_j^2=1\}.
\]
Then there is a tight frame for $\Reals^n$ consisting of $k$
vectors $\ub_1,\ldots,\ub_k\in\Ec$. \et
%%%%%%%%%%%%%%%%%%%%%%%%%%%%%%%%%%%

This result is valid for degenerate
ellipsoids (in which some of the major axes are infinitely
long). Our method of proof provides geometric insight to the problem, but
does not extend to infinite dimensions.

We note that, in the non-degenerate case, the definition of an ellipsoidal
 surface $\Ec$ given in Theorem \ref{construct},
is equivalent to the definition given in the introduction, specifying
 that the Hilbert space be $\Rn$.  Indeed, if $a_i > 1$
for all $i=1, \ldots, n$ and if $D = \txt{diag}(a_1, a_2, \ldots a_n)$, then
 $\sum_{i=1}^n a_ix_i^2=1$ iff $\ip{Dx}{x} = 1$ iff
$\|D^{\frac{1}{2}}x\|=1$ iff $D^{\frac{1}{2}}x \in \S(\Rone^n)$ iff $x \in
 D^{-\frac{1}{2}}\S(\Rn)$.  So $\Ec = \Ec_T$ for
$T = D^{-\frac{1}{2}}$, and thus $\Ec$ has the requisite form.  To reverse this
 argument for a non-diagonal positive operator $T$,
first diagonalize it by an orthogonal transformation given by rotations.
 Reversing the steps will then show that $\Ec_T$ is equivalent to
$\Ec$ for some choice of positive constants $\{a_1, \ldots a_n \}$.

The second theorem is used to prove Theorem \ref{exist} in the infinite dimensional case.  It has
 independent
interest to operator theory, and to our knowledge is a new result.  The proof, as well as a the
 corresponding
result in finite dimensions (Proposition
 \ref{rank1}), is contained in
Section \ref{operator}.  Some preliminaries are required before we state Theorem \ref{sum}.

It is well-known (see \cite{KR}) that a separably acting positive operator $A$ decomposes as the
 direct sum of a
positive operator $A_1$ with nonatomic spectral measure and a positive operator $A_2$ with purely
 atomic spectral
measure (i.e. a diagonalizable operator).  For $B \in \bh$,
 the
{\it essential norm} of $B$ is:
$$ \|B\|_{ess} := \inf \{\|B-K\|\;:\;  K \txt{is a compact operator in} \bh \}
 $$  In the proof of Proposition \ref{diag}, we have the special case where $A$ is a diagonal
 operator, $A = \txt{diag}(a_1, a_2,
 \ldots)$, with respect to some orthonormal
basis.  In this case, it is clear that $$\|A\|_{ess} = \sup\{\A > 0\;:\; |a_i| \geq \A \txt{for
 infinitely many}i \} $$

 For a positive operator $A$ with spectrum $\s(A)$, we have $\|A\| = \sup \{\L\,:\,\L \in \s(A) \}$
 and if $A$ is
 invertible, then $\|A^{-1}\|^{-1} = \inf\{\L\,:\, \L \in \s(A) \}$.  Similarly, $\|A\|_{ess}
 = \sup\{ \L\,:\, \L \in  \s_{ess}(A) \}$ and $\|A^{-1}\|^{-1}_{ess} = \inf \{\L\,:\, \L
\in \s_{ess}(A)\}$. In particular,  $\|A^{-1}\|^{-1}
 \leq \|A^{-1}\|^{-1}_{ess} \leq \|A\|_{ess} \leq \|A\|$.

For $A$ a positive operator, we say that $A$ has a {\it projection decomposition} if $A$ can be
expressed as the sum of a finite or infinite sequence of (not necessarily mutually orthogonal)
self-adjoint projections, with convergence in the strong operator topology.

%%%%%%%%%%%%%   Theorem 2    %%%%%%%%%%%%%%%%%%
 \bt\label{sum}  Let $A$ be a positive
operator in $\bh$ for $\hs$ a real or complex Hilbert
space with infinite dimension, and suppose
 $\|A\|_{ess}>1$.
Then $A$ has a projection decomposition.
\et
%%%%%%%%%%%%%%%%%%%%%%%%%%%%%%%%%%%%%%%%%%%%%%%%

Note that in this theorem, $A$ need not be invertible.  There are theorems in
 the literature (e.g. \cite{PT}) expressing operators as linear combinations
of projections and as sums of idempotents (non self-adjoint projections).  The
decomposition in Theorem \ref{sum} is different in that each term \textit{is} a
self-adjoint projection rather than a scalar multiple of a projection.

The next theorem states that every ellipsoidal surface contains a tight frame.
 We also include
some detailed information about the nature of the set of attainable frame
 bounds.

%%%%%%%%%  Theorem 3  %%%%%%%%%%%%%%%%%%
\bt\label{exist} Let T be a bounded invertible operator on a real or complex
Hilbert space.  Then the ellipsoidal surface $\Ec_T$ contains a tight frame. If
 $\hs$ is finite
dimensional with n = dim $\hs$, then for any integer $k \geq n$, $\Ec_T$ contains
a tight frame of length $k$, and every ETF on $\Ec_T$ of length $k$ has frame
bound $K= k \left[\tx{trace}(T^{-2})\right]^{-1}$. If $\txt{dim}\hs = \infty$, then for any constant
 $K >   \|T^{-2}\|^{-1}_{ess} $,
 $\Ec_T$ contains
a tight frame with frame bound $K$.  \et
%%%%%%%%%%%%%%%%%%%%%%%%%%%%%%%%%%%%%%%%%%

%%%%%%%%%%%%%%%%%%%%%%%%%%%%%%%%%%%%%%%%%%%%%%%%%%%%%%%%%%%%%%%%
%              Section                                         %
%%%%%%%%%%%%%%%%%%%%%%%%%%%%%%%%%%%%%%%%%%%%%%%%%%%%%%%%%%%%%%%%

\section{A Construction of ETFs in $\Rone^n$}\label{rotate}

We begin by showing that every ellipsoid can be scaled to contain an orthonormal
 basis.

\begin{lemma}\label{lem:onb}
Let $n\in\Nats$, let $a_1,\ldots,a_n\ge0$ be such that
$\sum_1^na_j=n$ and let
\[
\Ec=\{\xb=(x_1,\ldots,x_n)^t\in\Reals^n\mid\sum_1^na_jx_j^2=1\}.
\]
Then there is an orthonormal basis $\vb_1,\ldots,\vb_n$ for
$\Reals^n$ consisting of vectors $\vb_j\in\Ec$.
\end{lemma}
\begin{proof}
Proceed by induction on $n$. The case $n=1$ is trivial. Assume
$n\ge2$ and without loss of generality suppose $a_1\ge1$ and
$a_2\le1$. Let $\theta$ be such that
$a_1(\cos\theta)^2+a_2(\sin\theta)^2=1$ and let
$b_2=a_1(\sin\theta)^2+a_2(\cos\theta)^2$. Consider the rotation
matrix
\[
R=\left(\begin{matrix}
\cos\theta&\sin\theta \\
-\sin\theta&\cos\theta \\
&&1 \\
&&&\ddots \\
&&&&1
\end{matrix}\right).
\]
Then
\[
R^{-1}\Ec=\{(y_1,\ldots,y_n)^t\in\Reals^n\mid
y_1^2+2(a_1-a_2)y_1y_2\cos\theta\sin\theta+b_2y_2^2+\sum_3^na_jy_j^2=1\}.
\]
We have $b_2+\sum_3^na_j=n-1$. Let $\Vc$ be the subspace of
$\Reals^n$ consisting of all vectors of the form
$(0,x_2,\ldots,x_n)^t$ By the induction hypothesis, there is an
orthonormal basis $\ub_2,\ldots,\ub_n$ for $\Vc$ consisting of
vectors $\ub_j\in R^{-1}\Ec$. Let
$\ub_1=(1,0,\ldots,0)^t\in\Reals^n$, and let $\vb_j=R\ub_j$. Then
$\vb_1,\ldots,\vb_n$ is an orthonormal basis for $\Reals^n$
consisting of vectors $\vb_j\in\Ec$.
\end{proof}

In the case of a general ellipsoid, where $\sum_{j=1}^n a_j = r
>0$, the lemma gives a constant multiple of an orthonormal basis on
the ellipsoid.

\begin{proof}[Proof of Theorem \ref{construct}]
Consider the isometry $W:\Reals^n\to\Reals^k$ and the projection
$P=W^*:\Reals^k\to\Reals^n$ given by
\begin{align*}
W(x_1,\ldots,x_n)^t&=(x_1,\ldots,x_n,0,\ldots,0)^t \\
P(x_1,\ldots,x_k)^t&=(x_1,\ldots,x_n)^t.
\end{align*}
Let $a_j=0$ for $n+1\le j\le k$ and let
\[
\Ec'=\{\yb=(y_1,\ldots,y_k)^t\in\Reals^k\mid\sum_1^ka_jy_j^2=1\}.
\]

By Lemma~\ref{lem:onb}, there is a multiple of an orthonormal
basis $\vb_1,\ldots,\vb_k$ for $\Reals^k$ consisting of vectors
$\vb_j\in\Ec'$. Let $\ub_j=P\vb_j$. Then $\ub_j\in\Ec$. Moreover,
$\ub_1,\ldots,\ub_k$ is a tight frame for $\Reals^n$, because if
$\xb\in\Reals^n$, then
\[
\sum_{j=1}^k|\langle\xb,\ub_j\rangle|^2= \sum_{j=1}^k|\langle
W\xb,\vb_j\rangle|^2=\frac kr\|W\xb\|^2=\frac kr\|\xb\|^2.
\]
\end{proof}

\begin{remark}\label{otherproof} It is an elementary result in matrix theory \cite[Thm. 1.3.4]{HJ}
 that for any real $n \times n$
matrix $B$ acting on $\Rn$ there is an orthonormal basis $\{ u_1, \ldots u_n \}$
 for $\Rn$ so that the diagonal
elements $\ip{Bu_i}{u_i}$ of $B$ with respect to $\{u_1, \ldots u_n\}$ are all
 equal to $ \frac{1}{n} \left[ \tx{trace}(B) \right]$.
If we let $D = \txt{diag}(a_1, \ldots a_n )$, where the numbers $a_i$ are as in
 Lemma \ref{lem:onb}, then the condition $\ip{Dv}{v} = 1$
for a vector $v$ is exactly the condition for $v$ to be on the ellipsoid $\Ec$.
 Thus, letting $B=D$ and $v_i=u_i$ yields
another proof of Lemma \ref{lem:onb}.  The merit of the proof we give is that it
 is algorithmic and relates well to the paper.  It was
obtained by the second author in an undergraduate research (REU) program in
 which the other co-authors were mentors.
\end{remark}

%%%%%%%%%%%%%%%%%%%%%%%%%%%%%%%%%%%%%%%%%%%%%%%%%%%%%%%%%%%%%%%%
%              Section                                         %
%%%%%%%%%%%%%%%%%%%%%%%%%%%%%%%%%%%%%%%%%%%%%%%%%%%%%%%%%%%%%%%%

\section {Projection Decompositions for Positive Operators}\label{operator}

The arguments in the remainder of this paper hold for $\hs$ either a real or
 complex Hilbert space.

\begin{proposition} \label{rank1} Let $A \in \bh$
be a finite rank positive operator with integer trace $k$.  If
 $k \geq
rank(A)$, then $A$ is the sum of k projections of rank one.

 \end{proposition}

\bp We will construct unit vectors $x_1, x_2, \ldots x_k$ so that $A$ is the sum of the projections
 $x_i \otimes
x_i$ .  The proof uses induction on $k$. Let $n = \tx{rank}(A)$ and write $\hs_n = \tx{range}(A)$.
 If $k=1$, then
$A$ must itself be a rank-1 projection. Assume $k \geq 2$.  Select an orthonormal basis
 $\{e_i\}_{i=1}^n$ for
$\hs_n$ such that $A$ can be written on $\hs_n$ as a diagonal matrix with positive entries $a_1 \geq
 a_2 \cdots
\geq a_n
>0$.

Case 1: $(k > n)$ In this case, we have $a_1 > 1$, so we can take $x_k = e_1$. The
remainder on $\hs_n$ $$A- (x_k \otimes
x_k) = \tx{diag}(a_1-1, a_2, \ldots ,a_n)$$ has positive diagonal entries, still has rank $n$, and
now has  trace $k-1 \geq
n$. By the inductive hypothesis, the result holds.

Case 2:  $(k=n)$ We now have that $a_1 \geq 1$ and $a_n \leq 1$.  Given any finite rank,
self-adjoint $R \in \bh$,
let  $\mu_n(R)$ denote the $n$-th largest eigenvalue of $R$ counting multiplicity.
  Note that $\mu_n(A - (e_1 \otimes e_1)) \geq 0 $, $ \mu_n(A - (e_n \otimes e_n)) \leq 0$, and
$\mu_n(A - (x \otimes x))$ is a continuous function of $x \in \hs_n$.  Hence, there exists $y \in
 \hs_n$ such that
$\mu_n(A - (y \otimes y)) = 0$.  Choose $x_k = y$. Note the remainder $(A-(x_k \otimes x_k)) \geq 0$
 and
\begin{eqnarray*}
 \tx{trace}(A-(x_k \otimes x_k)) &=& n-1 \\ \tx{rank}(A-(x_k \otimes x_k)) &=& n
- 1 = k-1 \end{eqnarray*}
Again, by the inductive hypothesis, the result holds.

\ep

\begin{lemma}\label{rankk} Let $P_1, P_2, \ldots P_n$ be mutually
orthogonal projections on a Hilbert space $\hs$, all of
 the same nonzero
rank $k$, where $k$ can be finite or infinite.  Let $r_1, r_2, \ldots ,r_n$ be
 nonnegative
real numbers, and let $r = \sum_1^n r_i$.  Define the operator:
 $$ A= r_1P_1 + r_2P_2 + \cdots + r_nP_n$$
If the sum $r$ is an integer and $r
\geq n$, then there exist rank-$k$ projections  $Q_1, \ldots , Q_r$
such that: $$ A = Q_1 + Q_2 + \cdots + Q_r$$
\end{lemma}

 \bp
If $k=1$, then $r=\tx{trace}(A)$ and we have $\tx{rank}(A) \leq
n \leq r$, so the result follows from Proposition \ref{rank1}.  If $k
> 1$, each projection $P_i$ can be written as a sum of $k$
mutually orthogonal rank-1 projections: $$ P_i = P_{i1}+P_{i2}+
\cdots +P_{ik}$$ (Here and elsewhere in this proof, sums with indices running from $1$ to $k$
should be interpreted as infinite sums in the case where $k=\infty$.)  All rank-1 projections
$P_{ij}$ are thus mutually orthogonal.  Define operators $A_1, \ldots ,A_k$:  $$A_j=r_1P_{1j} +
r_2P_{2j}+\cdots +r_nP_{nj}$$ Now, $A = A_1 + \cdots + A_k$ and each $A_j$ has rank $n$ and trace
$r$.  By  Proposition \ref{rank1},
each $A_j$ can be written as a sum of $r$ rank-1 projections: $$A_j = T_{j1}+T_{j2}+\cdots +T_{jr}$$
Note that projections $T_{jl}$ and $T_{mp}$ are orthogonal if $j \neq m$.
 Define the rank-k projections
$Q_1, \ldots, Q_r$ by: $$Q_l = T_{1l}+T_{2l}+\cdots+T_{kl}$$
This gives $A = Q_1 + Q_2 + \cdots +Q_r$.
\ep

\bl\label{rational}  Let $A$ be a positive operator
 with finite spectrum contained in the rationals $\mathbb{Q}$, such that
all spectral projections are infinite dimensional, and also such that $\|A\|
 > 1$.  Then $A$ is a finite sum of self-adjoint projections.
\el

\bp  By hypothesis, there are mutually orthogonal infinite-rank projections $\,P_1, \ldots, P_n$ and
 positive
rational numbers $\,r_1 \geq r_2 \geq \cdots \geq r_n$ such that $$A= r_1P_1 + \cdots + r_nP_n$$
By hypothesis
$\|A\|>1$, hence $r_1 > 1$.

Write $r_i = s_i/t_i$ with $s_i$ and $t_i$ positive integers, and let $s=\sum_{i=1}^n s_i, \,
t=\sum_{i=1}^n t_i$.  We may assume $s \geq t$, for otherwise we can choose $m \in \Nats$ such that
$$ ms_1 + s_2 + \cdots + s_n \geq mt_1 + t_2 + \cdots + t_n$$
and replace $s_1$ with $ms_1$ and $t$ with $mt_1$.

Each $P_i$ can be written as a sum of $t_i$ mutually orthogonal infinite
 rank projections $P_{ij}\, ; \,j=1, \ldots
t_i$ which then allows us to write: $$ A =\sum_{i=1}^n \sum_{j=1}^{t_i} r_iP_{ij}$$
The operator is now a linear
combination of $\sum t_i = t$ mutually orthogonal
 projections of infinite rank, and the
 sum of the coefficients is now an integer $\sum t_i r_i = \sum s_i = s$.
 Since $s \geq t$,
 Lemma \ref{rankk} implies that A can be written as a sum of $s$ projections.

\ep

\begin{lemma}
Let $A$ be a positive operator which has a projection-decomposition. Then either $A$ is a
projection or  $\|A\|>1$.
\end{lemma}
\begin{proof}
Suppose, to obtain a contradiction, that $\|A\| \leq 1$ and that $A$ is not a projection.  By
 assumption, $A=\sum
P_i$ with the series converging strongly.  Thus $A-P_i \geq 0$ for all $i$.  Then $P_i(A-P_i)P_i
 \geq 0$, so
$P_iAP_i \geq P_i$.

Let $\mathcal{K}_i=P_i\hs$ and $B=P_iA|_{\mathcal{K}_i}$.  Then $B_i$ is positive and $B_i \geq
I_{\mathcal{K}_i}$ (the identity  operator on
$\mathcal{K}_i$).  Since $\|B_i\|\leq 1$, this implies $B_i = I_{\mathcal{K}_i}$, and thus $P_iAP_i
= P_i$.

Now, $P_i = P_i(\sum_j P_j)P_i = P_i + \sum_{j \neq i} P_iP_jP_i$, so $\sum_{j \neq i} P_iP_jP_i =
 0$.  Since each
$P_iP_jP_i \geq 0$, this implies $P_iP_jP_i = 0$.  Thus, $(P_jP_i)^*(P_jP_i) = 0$, so $P_jP_i=0$.
 Since this is
true for arbitrary $i, j$ with $i \neq j$, this shows that $A$ is the sum of mutually orthogonal
 projections, and
hence is itself a projection.  The contradiction proves the result.
\end{proof}

\begin{proposition}\label{nonatomic} Let $A$ be a positive operator in
$\bh$ with the property that all nonzero spectral projections for $A$ are of infinite rank.  If
 $\|A\|>1$, then
$A$ admits a projection decomposition as a sum of infinite rank projections.
\end{proposition}

\bp We will show that $A$ can be written as a sum $A=\sum_{i=1}^{\infty}A_i$ of
 positive operators, each satisfying the
hypotheses of Lemma \ref{rational}, where the sum converges in the strong
 operator topology.  We can then
decompose each of the operators $A_i$ as a finite sum of projections $A_{ij}$
 and then re-enumerate with a single index to
obtain a sequence $Q_i$ of projections which sum to $A$ in SOT.  Indeed, the
 partial sums of $\sum Q_i$ are dominated by $A$, hence
$\sum Q_i$ converges strongly to some operator $C$, and since the partial sums
 of $\sum A_i$ are also partial sums of $\sum Q_i$, the
sequence of partial sums of $\sum Q_i$ has a subsequence which converges to $A$,
 and hence $C=A$.

By hypothesis, we have $\|A\|>1$.  We may choose a positive rational number $\A > 1$ and a
nonzero spectral projection $G$ for $A$ such that $A \geq \A G$.  Let $B=A-\alpha G$, so that $B \geq 0$.  Using a standard argument, we can write
$B=\sum_{i=1}^{\infty} B_i$,   where each $B_i$ is a positive rational multiple of a spectral
projection for $A$, with convergence in the SOT.

We can write $G = \sum G_i$ as
 an infinite direct sum of nonzero infinite rank projections, with the requirement that $G_i$ be a subprojection of $G$ which commutes with all
the spectral projections for $A$. (This can clearly be done when the spectral projections for $A$ are all of
infinite rank.)  Now, let $A_i=\A G_i + B_i$.  We have $\|A_i\|\geq \A>1$.

By Lemma \ref{rational}, it follows that $A_i$ is a finite sum of projections.
 By the construction, we have the requisite form $A=\sum A_i$.

 \ep

\begin{proposition}\label{diag} Let $A$ be a positive operator in $\bh$ which is
diagonal with respect to some orthonormal basis $\{e_i\}$ for the
Hilbert space $\hs$.  Suppose  $\|A\|_{ess}>1$.  Then
 there is a sequence of rank-1
projections $\{P_i\}_{i=1}^{\infty}$ such that $A=\sum P_i$, where
the sum converges in the strong operator topology.
\end{proposition}

\bp Write $A$ as diag$(a_0, a_1, \ldots )$ and let $E_n=e_n \otimes e_n$.  Since $\|A\|_{ess}>1$,
 there is a
constant $\A > 1$ such that $a_i \geq \A $ for infinitely many i.  Let $k \geq 2$ be an integer such
that  $1+ \frac{2}{k-1} \leq \A$.  Permuting if necessary, we can without loss of generality assume
that the indices $n$ for which $a_n < \A$ are all multiples of $k$.

Let $B_0 = a_0E_0 + \cdots + a_{k-1}E_{k-1}$.
Therefore, we have rank$(B_0) \leq k$ and
\begin{eqnarray*}
\tx{trace}(B_0) &=& \sum_0^{k-1}a_i \\ &\geq&
 a_0 + (k-1)\A \\&\geq& a_0 + (k-1)\left(1+\frac{2}{k-1}\right)\\ &=&
 a_0+k+1
 \end{eqnarray*}
 Let $L_0$ be the greatest integer less than trace\,$(B_0)$.  Then
 $L_0 \geq k+1$.  Define $a'_{k-1}$ to be the real number $0 \leq
 a'_{k-1}\leq a_{k-1}$ such that if $$B'_0 = a_0E_0 + \cdots +
 a_{k-2}E_{k-2} + a'_{k-1}E_{k-1}$$ \noi then
 $$ \tx{trace}(B'_0) = L_0 \geq k+1 > \tx{rank}(B'_0)$$
By Proposition \ref{rank1}, $B'_0$ can be written as a sum of
$L_0$ rank-1 projections.

 In the next step, let $a''_{k-1} = a_{k-1}-a'_{k-1}$ and let
$$ B_1 = a''_{k-1}E_{k-1}+a_kE_k+a_{k+1}E_{k+1}+\cdots +
a_{2k-1}E_{2k-1}$$
Thus rank\,$(B_1) \leq k+1$ and
\begin{eqnarray*}
\tx{trace}(B_1) &=& a''_{k-1}+a_k+(a_{k+1}+ \cdots + a_{2k-1})\\
& \geq & a''_{k-1}+a_k+(k-1)\A \\
& \geq &a''_{k-1}+a_k+(k-1)\left( 1+ \frac{2}{k-1} \right) \\
&=& a''_{k-1}+a_k+k+1\\
&\geq& \tx{rank}(B_1)
\end{eqnarray*}

Construct $B_1'$ in a similar manner, so that its trace is an integer greater than or equal to its rank.  Then
$B'_1$ can be written as a sum of rank-1 projections using Proposition \ref{rank1}.

Proceeding recursively in a like manner, we may write $A=\sum_{j=1}^{\infty} B'_j$ converging in
SOT, where each $B'_j$ is a positive operator supported in $E_{jk-1}+ \cdots + E_{(j+1)k-1}$ and
with trace$(B'_j)$ an integer that is greater than or equal to rank$(B'_j)$.  Invoking Proposition
\ref{rank1} again to write each $B'_j$ as a sum of rank-1 projections, the proposition is proved.

\ep

\bp[Proof of Theorem \ref{sum}] Write $A=A_1 + A_2$, where $A_1$ and $A_2$ respectively denote the
 nonatomic and
purely atomic parts of $A$.  Then $\|A_1\|_{ess}=\|A_1\|$, and $\|A\|_{ess}=\max\{\|A_1\|, \|A_2\|_{ess} \}$.  So
$ \|A\|_{ess} > 1$ implies $ \|A_1\|>1 $ or $\|A_2\|_{ess} >1$.  Suppose first that $\|A_1\| > 1$.  Then there is
a nonzero spectral projection $P$ for $A_1$ and a constant
 $\A>1$ such that
$A_1 P \geq \A P$.  Let $Q$ be a nonzero spectral projection for $A_1$ dominated by $P$ such that $P-Q
 \neq 0$. Then
$A_1-\A Q$ satisfies the hypotheses of Proposition \ref{nonatomic}, so is projection decomposable.
 Also, $QA_2 =
A_2Q = 0$, so $A_2+\A Q$ is a diagonal operator with essential norm greater than or equal to $\A$,
 and so it is
projection decomposable by Proposition \ref{diag}.  The result follows by decomposing $A_1-\A Q$ and
 $A_2 + \A Q$
as sums of projections and combining the series.

For the case $\|A_1\| \leq 1$ and $ \|A_2\|_{ess} > 1$, we use a similar argument.  There is a constant
 $\A>1$
and an infinite rank spectral projection $P$ for $A_2$ such that $A_2-\A P \geq 0$.  Then $P$
 dominates a
projection $Q$ that commutes with $A_2$ such that both $Q$ and $P-Q$ are of infinite rank.  Then $A_2-\A Q$
 satisfies Proposition
\ref{diag} and hence has a projection decomposition.  The operator $A_1 + \A Q$ has norm greater
 than or equal to
$\A$ and all of its nonzero spectral projections have infinite rank, so it satisfies the hypotheses of Proposition
\ref{nonatomic}.  Thus, $A_1 + \A Q$ has a projection
 decomposition, and
we combine with the decomposition of $A_2-\A Q$ to get a projection decomposition for $A$. \ep

%%%%%%%%%%%%%%%%%%%%%%%%%%%%%%%%%%%%%%%%%%%%%%%%%%%%%%%%%%%%%%%%
%              Section                                         %
%%%%%%%%%%%%%%%%%%%%%%%%%%%%%%%%%%%%%%%%%%%%%%%%%%%%%%%%%%%%%%%%

\section{Ellipsoidal Tight Frames}\label{etf}

Let $\hs$ be a finite or countably infinite dimensional Hilbert space.  Let
$\{x_j\}_{j \in \J}$ be a frame for $\hs$, where $\J$ is some index set.
 Consider the standard
frame operator defined by: $$ Sw = \sum_{j \in \J} \ip{w}{x_j} x_j = \sum_{j \in
 \J} \left( x_j \otimes x_j \right)w $$
\noi Thus, $S = \sum_{\J} x_j \otimes x_j$, where this series of positive rank-1
 operators converges in the
strong operator topology (i.e. the topology of pointwise convergence).  In the
 special case where each
$\|x_j\|=1$, $S$ is the sum of the rank-1 projections $P_j = x_j \otimes x_j$.
 If we let $y_j = S^{-\frac 12}x_j$, then it is well-known
that $\{y_j\}_{j \in \J}$ is a Parseval frame (i.e. tight with frame bound 1).
 If each $\|x_j\|=1$, then $\{y_j\}_{j \in \J}$ is an
ellipsoidal tight frame for the ellipsoidal surface $\Ec_{S^{-\frac 12}}=S^{-\frac
 12}\S$.  Moreover, it is well-known (see \cite{HL}) that a sequence $\{x_j\}{j \in \J}
\subseteq \hs$ is a tight frame for $\hs$ if and only if the frame operator $S$
 is a positive scalar multiple of the identity,
i.e. $S=KI$, and in this case $K$ is the frame bound.

\begin{remark}  From the above paragraph, it is clear that a positive invertible operator is the
 frame operator for a frame
of unit vectors if and only if it admits a projection decomposition. (Each projection can be further decomposed
 into rank-1
projections, as needed.)
\end{remark}

The link between Theorem \ref{sum} and Theorem \ref{exist} is the following:

\begin{proposition}\label{proj} Let $T$ be a positive invertible operator in
 $\bh$, and let $K>0$ be a positive constant.  The ellipsoidal
surface $\Ec_T = T\S$ contains a tight frame $\{y_j\}$ with frame bound $K$ if
 and only if the operator $R = KT^{-2}$ admits a projection decomposition.
 In this case, $R$ is the frame
 operator for the spherical frame $\{T^{-1}y_j\}$.
\end{proposition}

 \begin{proof}
We present the proof in the infinite-dimensional setting, and note that the
 calculations in the finite dimensional
 case are identical but do not require discussion of convergence.  Let $\J$ be a
 finite or infinite index set.
 Assume $\Ec_T$ contains a tight frame $\{y_j\}_{j\in \J}$ with frame bound $K$.
   Then $\sum_{j\in \J} y_j \otimes y_j= KI$,
with the series converging in the strong operator topology.
   Let $x_j := T^{-1}y_j \in \S$, so $x_j \otimes x_j$ are projections.  We can
 then compute:
 \begin{eqnarray*}
 R = KT^{-2} &=& T^{-1} \left(
\sum_{j\in \J} y_j \otimes y_j \right) T^{-1}\\
&=& \sum_{j\in \J} T^{-1}y_j \otimes T^{-1}y_j
= \sum_{j\in \J} x_j \otimes x_j
\end{eqnarray*}

\noi This shows that $R$ can be decomposed as required.  Conversely, suppose $R$
 admits a projection decomposition $R = \sum P_j$, where $\{P_j\}$
are self-adjoint projections and convergence is in the strong operator topology.
  We can assume that the $P_j$ have rank-1, for otherwise we can decompose each
$P_j$ as a strongly convergent sum of rank-1 projections, and re-enumerate
 appropriately.  Since $P_j \geq 0$, the convergence is independent
of the enumeration used.  Write $P_j = x_j \otimes x_j$ for some unit vector
 $x_j$.  Letting $y_j = Tx_j$,
we have $y_j \in \Ec_T$, and we also have:
 \beqn
 KI = TRT &=& T \left( \sum_{j\in \J}
x_j \otimes x_j \right) T \\&=&\sum_{j\in \J}  Tx_j \otimes Tx_j
= \sum_{j\in \J}   y_j \otimes y_j
 \eeqn
This shows that $\sum y_j \otimes y_j$ converges in the strong operator topology
 to $KI$.  Thus, $\{y_j\}_{j\in \J}$ is a
tight frame on $\Ec_T$, as required.

\end{proof}

 \bp[Proof of Theorem \ref{exist}]  Let $\Ec$ be an ellipsoid.  Then $\Ec =
 \Ec_T = T\S$ for some positive invertible $T \in \bh$.
 Let $K$ be a positive constant, and let $R = KT^{-2}$.

The condition  $K>\|T^{-2}\|_{ess}^{-1}$ implies
 $\|R\|_{ess}>1$.  So, by Theorem \ref{sum}, $R$ admits a projection
 decomposition, and thus Proposition \ref{proj} implies
 that $\Ec$ contains a tight frame with frame bound $K$.

In the finite dimensional case, let $n=\txt{dim}\hs$.  Proposition \ref{proj}
 states that $\Ec$ will contain a tight frame with
frame bound $D$ if and only if $KT^{-2}$ admits a projection decomposition, and
 by Proposition \ref{rank1} this happens if and only if
 trace\,$(KT^{-2})$ is an integer $k \geq n$, and in this case $k$ is the length
 of the frame.  Thus, we have $K = k[\tx{trace}(T^{-2})]^{-1}$.
 Therefore, every ellipsoid $\Ec = \Ec_T$ contains a tight frame of every length
 $k \geq n$, and every such tight frame has
 frame bound $k[\tx{trace}(T^{-2})]^{-1}$.
 \ep

\begin{corollary}
Every positive invertible operator $S$ on a separable Hilbert space $\hs$ is the frame operator for a spherical
frame. If $\hs$ has finite dimension $n$, then for every integer $k \geq n$, $S$ is the frame operator for a spherical
frame of length $k$, and the radius of the sphere is $ \sqrt{ \frac{\tx{trace}(S)}{k}}$. If $\hs$ is
infinite-dimensional, the radius of the sphere can be taken to be any positive number $r<\|S\|_{\tx{ess}}^{\frac
12}$.
\end{corollary}

\begin{proof}
In the finite dimensional case, let $c = \frac{k}{\tx{trace}(S)}$ and $A= cS$, so that $\tx{trace}(A) = k$. Then,
by Proposition \ref{rank1}, $A$ has a projection decomposition into $k$ rank-1 projections, making $A$ the frame
operator for the frame of unit vectors $\{x_i\}_{i=1}^k$.  Thus, $S$ is the frame operator for $ \left\{
\frac{x_i}{\sqrt{c}} \right\}_{i=1}^k$.

When $\hs$ has infinite dimension, let $c$ be any constant greater than $\|S\|_{\tx{ess}}^{-1}$, and let $A=cS$.
The hypotheses
 of Theorem
\ref{sum} are satisfied, so $A$ admits a projection decomposition.  Then $A$ is the frame operator
 for a frame
$\{x_i\}$ of unit vectors, so $S$ is the frame operator for the spherical frame $ \displaystyle
 \left\{
\frac{x_i}{\sqrt{c}} \right\}$.
\end{proof}

\begin{remark}\label{others}  We know of at least two groups who have independently and
simultaneously  proved our finite-dimensional ellipsoidal tight frame results.
 Paulsen and Holmes have a proof similar to the discussion in Remark \ref{otherproof}. \cite{HP}
 Casazza and Leon have shown in \cite{CL} the existence of
 ``spherical frames for $\Rn$ with a given frame operator'', which is an
 equivalent problem.

 \end{remark}

%%%%%%%%%%%%%%%%%%%%%%%%%%%%%%%%%%%%%%%%%%%%%%%%%%%%%%%%%%
%     Bibliography                                       %
%%%%%%%%%%%%%%%%%%%%%%%%%%%%%%%%%%%%%%%%%%%%%%%%%%%%%%%%%%
\providecommand{\bysame}{\leavevmode\hbox to3em{\hrulefill}\thinspace}


\begin{thebibliography}{10}

\bibitem{BF}
John~J. Benedetto and Matthew Fickus, \emph{Finite normalized tight frames},
  Preprint, 2001.

\bibitem{CKLT}
Peter~G. Casazza, Jelena Kova\v{c}evi\'{c}, Manuel Leon, and Janet~C. Tremain,
  \emph{Custom built tight frames}, Preprint, 2002.

\bibitem{CL}
Peter~G. Casazza and Manuel~T. Leon, \emph{Frames with a given frame operator},
  Preprint, 2002.

\bibitem{D}
Ingrid Daubechies, \emph{Ten lectures on wavelets}, CBMS-NSF Regional
  Conference Series in Applied Mathematics, Society for Industrial and Applied
  Mathematics (SIAM), Philadelphia, PA, 1992.

\bibitem{DGM}
Ingrid Daubechies, A.~Grossmann, and Y.~Meyer, \emph{Painless nonorthogonal
  expansions}, J. Math. Phys. \textbf{27} (1986), no.~5, 1271--1283.

\bibitem{DS}
R.J. Duffin and A.C. Schaeffer, \emph{A class of nonharmonic {Fourier} series},
  Trans. Amer. Math. Soc. \textbf{72} (1952), no.~2, 341--366.

\bibitem{GKK}
Vivek~K. Goyal, Jonathan~A. Kelner, and Jelena Kova\v{c}evi\'{c},
  \emph{Quantized frame expansions with erasures}, Appl. Comput. Harmon. Anal.
  \textbf{10} (2001), 203--233.

\bibitem{HL}
Deguang Han and David~R. Larson, \emph{Frames, bases, and group
  representations}, Mem. Amer. Math Soc. \textbf{147} (2000), no.~697, 1--94.

\bibitem{HW}
Eugenio Hern\'{a}ndez and Guido Weiss, \emph{A first course in wavelets}, CRC
  Press, LLC, Boca Raton, FL, 1996.

\bibitem{HP}
Roderick~B. Holmes and Vern~I. Paulsen, \emph{Optimal frames for erasures},
  Preprint, 2003.

\bibitem{HJ}
R.~Horn and C.~Johnson, \emph{Topics in matrix analysis}, Cambridge University
  Press, 1991.

\bibitem{KR}
Richard~V. Kadison and John~R. Ringrose, \emph{Fundamentals of the theory of
  operator algebras, {V}olume {I}: Elementary theory}, Graduate Studies in
  Mathematics, American Mathematical Society, Providence, RI, 1997.

\bibitem{PT}
Carl Pearcy and David Topping, \emph{Sums of small numbers of idempotents},
  Michigan Math. J. \textbf{14} (1967), 453--465.

\end{thebibliography}
\end{document}